\input amstex\documentstyle{amsppt}  
\pagewidth{12.5cm}\pageheight{19cm}\magnification\magstep1
\topmatter
\title Discretization of Springer fibres\endtitle
\author G. Lusztig\endauthor
\address{Department of Mathematics, M.I.T., Cambridge, MA 02139}\endaddress
\thanks{Supported by NSF grant DMS-1855773.}\endthanks
\endtopmatter   
\document

\define\bF{\bar F}

\define\frl{\forall}

\define\si{\sim}

\define\sqc{\sqcup}

\define\qua{\quad}

\define\op{\oplus}
   
\define\part{\partial}
\define\emp{\emptyset}

\define\n{\notin}

\define\m{\mapsto}
\define\do{\dots}

\define\sm{\smallmatrix}
\define\esm{\endsmallmatrix}
\define\sub{\subset}    

\define\T{\times}
\define\ti{\tilde}
\define\nl{\newline}
\redefine\i{^{-1}}

\define\un{\underline}

\define\ot{\otimes}

\define\Ad{\text{\rm Ad}}

\redefine\b{\beta}

\redefine\d{\delta}
\define\e{\epsilon}

\define\io{\iota}
\redefine\o{\omega}

\define\ph{\phi}

\define\r{\rho}
\define\s{\sigma}
\redefine\t{\tau}

\redefine\l{\lambda}
\define\z{\zeta}
\define\x{\xi}

\redefine\G{\Gamma}

\define\Om{\Omega}

\define\BB{\bold B}
\define\CC{\bold C}

\define\NN{\bold N}

\define\QQ{\bold Q}

\define\ZZ{\bold Z}

\define\ca{\Cal A}
\define\cb{\Cal B}

\define\cf{\Cal F}
\define\cg{\Cal G}

\define\cm{\Cal M}

\define\cz{\Cal Z}

\define\fb{\frak b}

\define\fg{\frak g}

\define\fp{\frak p}

\define\fC{\frak C}

\define\tb{\ti b}

\head Introduction\endhead
\subhead 0.1\endsubhead
Let $G$ be an almost simple simply connected algebraic group over $\CC$ with Lie
algebra $\fg$. 
Let $e\in\fg$ be a fixed nilpotent element and let $\cb_e$ be the variety of Borel
subalgebras of $\fg$ that contain $e$ (a Springer fibre). We fix a homomorphism of algebraic groups
$\z:SL_2(\CC)@>>>G$ whose differential carries $\left(\sm 0&1\\0&0\esm\right)$ to $e$. 
Let $F$ be the centralizer in $G$ of the image of $\z$ (a reductive group). Let $\bF=F/(F^0\cz_G)$.
(For any algebraic group 
$\cg$ we denote by $\cg^0$ the identity component of $\cg$; $\cz_G$ is the centre of $G$.)
Following \cite{L83} we view $\cb_e$ as a variety with $\CC^*$-action given by
$\l:\fb\m\Ad(\z\left(\sm\l&0\\0&\l\i\esm\right))\fb$.

Let $W$ be the (extended) affine Weyl group corresponding to the dual of $G$.
Let $c$ be the two-sided cell of $W$ associated to
the $G$--conjugacy class of $u=\exp(e)\in G$ in \cite{L89, 4.8}.
In this paper we consider the following four sets associated to $e$.

(a) The subset $\un\BB_{\cb_e}^\pm$ of $K_{\CC^*}(\cb_e)$ (the $K$-group of $\CC^*$-equivariant 
coherent sheaves on $\cb_e$) introduced in \cite{L99a, 5.15}.

(b) The set $R(c)$ of right cells of $W$ that are contained in $c$.

(c) The set $\Xi_e$ of connected components of the fixed point set $\cb_e^{\CC^*}$ of the
$\CC^*$-action on $\cb_e$.

(d) The set $\bar\Xi_e$ of orbits of the $\bF$-action on $\Xi_e$ induced by the conjugation
action of $F$ on $\cb_e^{\CC^*}$.
\nl
In the rest of this paper $\un\BB_{\cb_e}^\pm$ is renamed as $\BB^\pm_e$. 
Note that in \cite{L99} it is conjectured (and in \cite{BM} it is proved) that

(e) $\BB_e^\pm$ is a signed basis of the $K_{\CC^*}(point)$-module $K_{\CC^*}(\cb_e)$.
\nl
One of the themes of this paper is a conjectural diagram involving the sets (a)-(d).
$$\CD
\BB_e @>\r>> R(c)\\
@V\s VV          @V\s'VV\\
\Xi_e@>\r'>>\bar\Xi_e
\endCD$$
Here $\BB_e$ is the set of orbits of multiplication by $\{1,-1\}$ on $\BB_e^\pm$;

(f) $\r$ is the (conjectural) map in \cite{L99, 17.1(c)} which identifies $R(c)$ with the set
of $\bF$-orbits on $\BB_e$ (for the action of $\bF$ on $\BB_e$ induced by the conjugation
action of $F$ on $\cb_e$);
\nl
$\s$ is a (conjectural)
surjective map (compatible with the actions of $\bF$) discussed in Section 1; $\r'$ is the obvious
orbit map; $\s'$ is the unique (surjective) map which makes the diagram commutative.

In this paper we introduce a new (conjectural) signed basis $\ti\BB^\pm_e$ of (a localization of)
$K_{\CC^*}(\cb_e)$ which is in natural bijection with $\BB^\pm_e$ and is such that
$\BB^\pm_e$ can be reconstructed from the knowledge of
$\ti\BB^\pm_e$ and from the bar-involution of $K_{\CC^*}(\cb_e)$
in a way similar (but more intricate) to the way the canonical basis of the $+$ part of
a quantum group can be reconstructed from a PBW basis of that $+$ part. Thus we can think of
$\ti\BB^\pm_e$ as being something like a PBW (signed) basis. The set 
$\ti\BB^\pm_e$ is naturally partitioned into subsets indexed by $\Xi_e$ in (c); this can
be viewed as a surjective map $\BB^\pm_e@>>>\Xi_e$ which factors through a
surjective map $\BB_e@>\s>>\Xi_e$ appearing in the diagram above.

\subhead 0.2\endsubhead
The set $\BB_e$ is a {\it discretization} (or discrete analogue) of $\cb_e$ in the sense
that it is a finite set with a number of elements equal to the sum of Betti numbers
(or equivalently the sum of Betti numbers in even degrees) of $\cb_e$. (This follows from 0.1(e).)

\subhead 0.3\endsubhead
The set $\BB_e$ appears in representation theory in at least two different
ways. It indexes the simple objects in a certain block of unrestricted
representations of the analogue of $\fg$ over a field of positive, large characteristic 
(this has been conjectured in \cite{L98, \S14} and proved in \cite{BM}).

A second application of the set $\BB_e$ is as follows.
Now $F$ acts on $\BB_e$ via its quotient $\bF$.
By \cite{L99a, 17.1}, the $F$-set $\BB_e$ should be
the $F$-set $Y$ appearing in the conjecture \cite{L89, 10.5} which provides

(a) a bijection between $c$ and the set of indecomposable 
$F$-equivariant vector bundles on $Y\T Y$ up to isomorphism.
\nl
For $G$ of type $A$, (a) has been established in \cite{XI}; for general $G$, a weak
form of (a) has been established in \cite{BFO}.

\subhead 0.4\endsubhead
In section 2 we state some conjectures which, if true, would describe completely
the finite set $\BB_e$ with action of $\bF$ (that is, they describes which isotropy groups
appear and how many points have isotropy groups in a fixed conjugacy class).
In section 3 we illustrate the bijection 0.3(a) in an example.

\head 1. The maps $\BB_e@>>>\Xi_e$, $R(c)@>>>\bar\Xi_e$  \endhead
\subhead 1.1\endsubhead
Let $\cb$ be the variety of Borel subalgebras of $\fg$.
We have $\cb_e=\{\fb\in\cb;e\in\fb\}$. As in 0.1 we consider
$K_{\CC^*}(\cb_e)$, the $K$-theory of $\CC^*$-equivariant coherent
sheaves on $\cb_e$; we denote it by $K_e$. We regard $K_e$ as a module over $\ca:=\ZZ[v,v\i]$
(the representation  ring of $\CC^*$) in the usual way.
Here $v$ is an indeterminate representing the identity homomorphism $\CC^*@>>>\CC^*$.

In \cite{L99a, 5.15} we have defined an involution  $\ti\b:K_e@>>>K_e$, a
symmetric $\ca$-bilinear pairing $(||):K_e\T K_e@>>>\ca$ and the subset
$$\un\BB_{\cb_e}^\pm=\{\x\in K_e);\ti\b(\x)=\x,(\x||\x)\in 1+v\i\ZZ[v\i]\}$$
of $K_e-\{0\}$ (now denoted by $\BB_e^\pm$). We will also write $\bar{}$ instead of $\ti\b$.)

\subhead 1.2\endsubhead
Let ${}'\ca$ be the subring of $\QQ(v)$ consisting of quotients
$f/g$ where $g\in\ZZ[v]$ has constant term $1$ and $f\in\ca$; 
let ${}''\ca$ be the subring of $\QQ(v)$ consisting of quotients
$f/g$ where $g\in\ZZ[v]$ has constant term $1$ and $f\in\ZZ[v]$.
We have $\ca\sub{}'\ca,{}''\ca\sub{}'\ca$.

For any $m\in\ZZ$ let 
$$\fg_m=\{x\in\fg;\Ad(\z\left(\sm\l&0\\0&\l\i\esm\right))x=\l^mx\qua\frl\l\in\CC^*\}.$$
Then $\fp:=\sum_{m\in\NN}\fg_m$ is the Lie algebra of a parabolic subgroup
$P$ of $G$ containing $F$. Let $\cm$ be the (finite) set of orbits of $P$ on $\cb$ (for the
conjugation action). Let $\cm_e=\{\o\in\cm;\cb_e\cap\o\ne\emp\}$.
Let $\ti\cm_e$ be the set of all subvarieties $X\sub\cb_e$ such that $X$ is
a connected component of $\cb_e\cap\o$ for some $\o\in\cm_e$.

If $X\in\ti\cm_e$ is a connected component of $\cb_e\cap\o$ with $\o\in\cm_e$,
we set $\cb_e^{<X}=\cup_{\o'}(\cb_e\cap\o')$
where $\o'\in\cm_e$ is subject to $\o'\sub\bar\o$ (closure in $\cb$) and $\o'\ne\o$;
we set $\cb_e^{\le X}=X\cup\cb_e^{<X}$. 

By arguments in \cite{L99(a), Section1} (based on results in \cite{DLP}) we see that
the $\ca$-linear maps $K_{\CC^*}(\cb_e^{<X})@>>>K_e$,
$K_{\CC^*}(\cb_e^{\le X})@>>>K_e$ induced by the closed imbedding $\cb_e^{<X}\sub\cb_e$,
$\cb_e^{\le X}\sub\cb_e$, are injective; hence the ${}'\ca$-linear maps
${}'\ca\ot_\ca K_{\CC^*}(\cb_e^{<X})@>>>{}'\ca\ot_\ca K_e$,
${}'\ca\ot_\ca K_{\CC^*}(\cb_e^{\le X})@>>>{}'\ca\ot_\ca K_e$
obtained by extension of scalars are injective.
Hence ${}'\ca\ot_\ca K_{\CC^*}(\cb_e^{<X})$ and ${}'\ca\ot_\ca K_{\CC^*}(\cb_e^{\le X})$
can be identified with their image ${}'K^{<X}_e$ and ${}'K^{\le X}_e$ in
${}'K_e:={}'\ca\ot_\ca K_e$.) The same arguments show that we have an exact sequence
The same arguments show that we have an exact sequence
$$0@>>>K_{\CC^*}(\cb_e^{<X})@>>>K_{\CC^*}(\cb_e^{\le X})@>>>K_{\CC^*}(X)@>>>0$$
associated to the inclusions $\cb_e^{<X}\sub\cb_e^{\le X}$, $X\sub\cb_e^{\le X}$;
from this we deduce an exact sequence $0@>>>{}'K^{<X}_e@>>>{}'K^{\le X}_e@>t>>{}'K_{\CC^*}(X)$
where ${}'K_{\CC^*}(X)={}'\ca\ot_\ca K_{\CC^*}(X)$.

We have naturally $K_e\sub{}'K_e$.

Now $(||):K_e\T K_e@>>>\ca$ extends to a symmetric ${}'\ca$-bilinear pairing
${}'K_e\T{}'K_e@>>>{}'\ca$. For any $X\in\ti\cm_e$ let 
${}'K^X_e$ be the set of all $\x\in{}'K^{\le X}_e$ such that
$(\x||\x')=0$ for any $\x'\in{}'K^{<X}_e$.
Restricting $t:{}'K^{\le X}_e@>t>>{}'K_{\CC^*}(X)$ to ${}'K^X_e$ we obtain a map

(a) ${}'K^X_e@>>>{}'K_{\CC^*}(X)$.
\nl
Let ${}''K_e$ be the ${}''\ca$-submodule of ${}'K_e$ generated by $\BB^{\pm}_e$.

\subhead 1.3 \endsubhead
We now state some conjectural properties of the submodules ${}'K^X_e$ of ${}'K_e$.

(i) We have a direct sum decomposition ${}'K_e=\op_{X\in\ti\cm_e}{}'K^X_e$. Hence for
$b\in\BB_e^{\pm}$ we can write uniquely $b=\sum_{X\in\ti\cm_e}b^X$ where $b^X\in{}'K^X_e$.
Moreover, the maps 1.2(a) are isomorphisms, hence they convert the direct sum
decomposition above into ${}'K_e=\op_{X\in\ti\cm_e}{}'K_{\CC^*}(X)$. 

(ii) Let $b\in\BB_e^{\pm}$. There is a
unique $X_b\in\cm_e$ such that $b^X\in v({}''K_e)$
for all $X\in\ti\cm_e-\{X_b\}$ and $b^{X_b}-b\in v({}''K_e)$ (so that
$b^{X_b}\n v({}''K_e)$). The
map $\BB^{\pm}_e@>>>\ti\cm_e$, $b\m X_b$ is surjective.
\nl
In this and the next subsection (but not in other subsections)
we identify $\BB_e$ with a subset of $\BB_e^\pm$ by choosing one element in each
orbit of multiplication by $\{1,-1\}$ on $\BB_\e^\pm$.
Setting $\ti b=b^{X_b}$ for any $b\in\BB_e$ we have
$\ti b=\sum_{b'\in\BB_e}c_{b,b'}b'$ where $c_{b,b'}\in{}''\ca$ satisfy
$c_{b,b}\in 1+v({}''\ca)$, $c_{b,b'}\in v({}''\ca)$ for $b\ne b'$. It follows that the square
matrix $(c_{b,b'})$ indexed by $\BB_e\T\BB_e$ has determinant in $1+v({}''\ca)$
hence is
invertible in ${}''\ca$. Since $\{b;b\in\BB_e\}$ is an ${}''\ca$ basis of ${}''K_e$,
it follows that

(a) $\ti\BB_e:=\{\ti b;b\in\BB_e\}$ is again an ${}''\ca$ basis of ${}''K_e$.

\subhead 1.4\endsubhead
We show that the $\ca$-basis $\BB_e$ van be reconstructed from
the ${}'\ca$-basis $\ti\BB_e$ of ${}'K_e$ (assuming 1.3(i),(ii)).

We shall indicate a number of steps which start with $\ti\BB_e$ and end with $\BB_e$
(the definition of these steps does not involve $\BB_e$, but the verification of their
correctness does).

Step 1. We note that ${}''K_e$ is defined purely in terms of $\ti\BB_e$ (it is
the ${}''\ca$-submodule of ${}'K_e$ generated by $\ti\BB_e$).

Step 2. We set ${}^+K_e=K_e\cap{}''K_e$.

Step 3. Let ${}^-K_e$ be the image of ${}^+K_e$ under $\bar{}:K_e@>>>K_e$.

Step 4. We form ${}^+K_e\cap{}^-K_e$.

Step 5. We have a map ${}^+K_e\cap{}^-K_e@>\io>>{}^+K_e/v{}^+K_e$ (restriction of the
obvious map ${}^+K_e@>>>{}^+K_e/v{}^+K_e$).

Step 6. We have a map ${}^+K_e/v{}^+K_e@>\io'>>{}''K_e/v{}''K_e$ induced by
the obvious inclusion ${}^+K_e\sub{}''K_e$.

Step 7. For any $\fb\in\ti\BB_e$ there is a unique element $\t(\fb)\in
{}^+K_e\cap{}^-K_e$
such that $\io'\io(\t(\fb))$ is the image of $\fb$ in ${}''K_e/v{}''K_e$.

Step 8. The elements $\{t(\fb);\fb\in\ti\BB_e\}$ form a $\ZZ$-basis of 
${}^+K_e\cap{}^-K_e$ and an $\ca$-basis of $K_e$. This is $\BB_e$.
\nl
We now justify Step 7. Note that ${}^+K_e$ is the set of all
$\sum_{b\in\BB_e}c_bb$ where for any $b$ we have
$c_b\in\ca\cap{}''\ca$ or equivalently $c_b\in\ZZ[v]$.
(If $a\in\ZZ[v,v\i]$ is of the form $f/g$ where $g\in\ZZ[v]$ has constant term $1$ and
$f\in\ZZ[v]$, then $a\in\ZZ[v]$. Indeed, we have $a=\sum_{i\in\ZZ}a_iv^i$ where $a_i\in\ZZ$
satisfies $a_i=0$ for $i\gg0$ and for $i\ll0$, since $a\in\ca$, and $a_i=0$ for $i<0$, since
$a\in{}''\ca$.) It follows that ${}^-K_e$ is the set of all $\sum_{b\in\BB_e}c_bb$ where
for any $b$ we have $c_b\in\ZZ[v\i]$. Hence ${}^+K_e\cap{}^-K_e$  
is the set of all $\sum_{b\in\BB_e}c_bb$ where for any $b$ we have $c_b\in\ZZ$.
The map $\io'$ in Step 6 is an isomorphism. 
(We use that the map $\ZZ[v]/v\ZZ[v]@>>>{}''\ca/v({}''\ca)$ induced by the inclusion
$\ZZ[v]@>>>{}''\ca$ is an isomorphism.) Moreover the map $\io$ in Step 5 is an isomorphism.
Now Step 7 holds in view of Steps 5 and 6.
We now justify Step 8. Under the isomorphism $\io'\io$, the $\ZZ$-basis $\BB_e$ of 
${}^+K_e\cap{}^-K_e$ corresponds to the $\ZZ$-basis of ${}''K_e/v({}''K_e)$
formed by the image of $\BB_e$ or equivalently by the image of $\ti\BB_e$.
This justifies Step 8.

We note that we can reconstruct $\BB_e$ from slightly less than the knowledge of $\ti\BB_e$: it is
enough to have ${}''K_e$ and the image of $\ti\BB_e$ under
${}''K_e@>>>{}''K_e/v({}''K_e)$.

\subhead 1.5\endsubhead
In this subsection we assume that $G$ is of type $A_2$ and $e\in\fg$ is subregular
nilpotent. Using \cite{L02, Sec.5}, we see that
$\BB_e^{\pm}$ consists of $\pm$ three elements $b_1,b_2,b_3$
satisfying $(b_i||b_i)=1+v^{-2}$ and $(b_i||b_j)=-v\i$ if $i\ne j$.
The set $\ti\cm_e$ has three elements which can be denoted by $X_1,X_2,X_3$ so that
${}'K^{X_3}_e={}'K^{\le X_3}_e$ has basis $\{b_3+vb_1+vb_2\}$,
${}'K^{\le X_i}_e$ has basis $\{b_3+vb_1+vb_2,b_i\}$ for $i=1,2$. It follows that
for $i=1,2$, ${}'K^{X_i}_e$ has basis $\{b_i-\d\i(v^3-v^2)(b_3+vb_1+vb_2)\}$
where $\d=1-v^2-2v^3+2v^4$.

We have
$$b_i=\d\i(v^3-v^2)(b_3+vb_1+vb_2)+(b_i-\d\i(v^3-v^2)(b_3+vb_1+vb_2)),$$
for $i=1,2$,
$$\align&b_3=\d\i(1-v^2)(b_3+vb_1+vb_2)-v(b_1-\d\i(v^3-v^2)(b_3+vb_1+vb_2))\\&
-v(b_2-\d\i(v^3-v^2)(b_3+vb_1+vb_2)).\endalign$$
Hence
$$\tb_i=b_i-\d\i(v^3-v^2)(b_3+vb_1+vb_2)\text { for }i=1,2,$$
$$\tb_3=\d\i(1-v^2)(b_3+vb_1+vb_2).$$
\nl
The map $\BB_e^{\pm}@>>>\ti\cm_e$ is $\pm b_i\m X_i$ for $i=1,2,3$.
We see that 1.3(i),(ii) hold in this case.

\subhead 1.6\endsubhead
In this subsection we assume that $G$ is of type $D_4$ or $G_2$ and $e\in\fg$ is subregular
nilpotent. Using \cite{L99b}, \cite{L02}, we see that $\BB_e^{\pm}$ consists of $\pm$
five elements
$b_0,b_1,b_2,b_3,b_4$ satisfying $(b_i||b_i)=1+v^{-2}$ for $i=0,1,2,3,4$,
$(b_i||b_j)=0$ if $i\ne j$ in $1,2,3,4$, $(b_0||b_i)=-v\i$ for $i=1,2,3,4$.
The set $\ti\cm_e$ has four elements which can be denoted by $X_0,X_1,X_2,X_3$ so that
${}'K^{X_0}_e={}'K^{\le X_0}_e$ has basis $\{b_0,b_4+v^2(b_1+b_2+b_3)\}$.
${}'K^{\le X_i}_e$ has basis $\{b_4+v^2(b_1+b_2+b_3),b_0,b_i\}$ for $i=1,2,3$.
It follows that for $i=1,2,3$, ${}'K^{X_i}_e$ has basis 
$$\{b_i+(v+v^3)\e\i b_0+v^4\e\i(b_4+v^2(b_1+b_2+b_3))\}$$
where $\e=1+2v^2-3v^6$. We have
$$\align&b_4=\e\i(1+2v^2)(b_4+v^2(b_1+b_2+b_3))+3\e\i(v^3+v^5)b_0\\&
-\sum_{i\in\{1,2,3\}}v^2\e\i(v^4(b_4+v^2(b_1+b_2+b_3))+(v+v^3)b_0+\e b_i),\endalign$$
$$\align&b_i=-\e\i v^4(b_4+v^2(b_1+b_2+b_3))-\e\i(v+v^3)b_0\\&
+\e\i(v^4(b_4+v^2(b_1+b_2+b_3))+(v+v^3)b_0+\e b_i)\endalign$$
for $i=1,2,3$. Hence
$$\tb_0=b_0,$$
$$\tb_4=\e\i(1+2v^2)(b_4+v^2(b_1+b_2+b_3))+3\e\i(v^3+v^5)b_0,$$
$$\tb_i=\e\i(v^4(b_4+v^2(b_1+b_2+b_3))+(v+v^3)b_0+\e b_i)\text{ for }i=1,2,3.$$
The map $\BB_e^{\pm}@>>>\ti\cm_e$ is $\pm b_i\m X_i$ for $i=0,1,2,3$ and $\pm b_4\m X_0$.
We see that 1.3(i),(ii) hold in this case.

\subhead 1.7\endsubhead
If $\o\in\cm_e$, then $F$ acts on $\cb_e\cap\o$ by conjugation. This induces an action of
$\bF$ on the set of connected components of $\cb_e\cap\o$. By \cite{DLP}, this action of $\bF$ is
transitive. Thus, $\bF$ acts naturally on $\ti\cm_e$ and the map $\ti\cm_e@>>>\cm_e$ (with
$X\m\o$ when $X\sub\cb_e\cap\o$) has fibres given precisely by the $\bF$-orbits on $\ti\cm_e$.

By \cite{DLP}, if $X\in\ti\cm^e$, then $X^{\CC^*}=X\cap\cb_e^{\CC^*}$ is a connected
component $\cb_e^{\CC^*}$ that is an element of $\Xi_e$; moreover, $X\m X^{\CC^*}$
is a bijection $\ti\cm_e@>\si>>\Xi_e$. Thus we may identify $\ti\cm_e$ with $\Xi_e$
and $\cm_e$ with $\bar\Xi_e$ (see 0.1).

Using the identification $\ti\cm_e=\Xi_e$, the map $\BB^{\pm}_e@>>>\ti\cm_e$ in 1.3(ii) can be
identified with a map $\BB^{\pm}_e@>>>\Xi_e$, which factors through a (surjective) map
$\s:\BB_e@>>>\Xi_e$. Thus all maps in the diagram in 0.1 are defined.

\subhead 1.8\endsubhead
One can define a (non-conjectural)
direct sum decomposition $\QQ(v)\ot_\ca K_e=\op_{X\in\ti\cm_e}K(X)$
into $\QQ(v)$-vector subspaces $K(X)$
indexed by $X\in \ti\cm_e$ by noting that by a known localization
property we have $\QQ(v)\ot_\ca K_e=\QQ(v)\ot_\ca K_{\CC^*}(\cb_e^{\CC^*})$
and then using the direct sum decomposition of the last vector space coming from the
decomposition of $\cb_e^{\CC^*}$ into connected components (which are indexed by $\ti\cm_e$).
One can project any $b\in\BB^\pm_e$ to the summands in this decomposition and one can ask whether
these projections behave as in 1.3(ii). It appears that this is not the case.

\head 2. $\BB_e$ and the Burnside group of $\bF$\endhead
\subhead 2.1\endsubhead
Let $H$ be a finite group.
Let $\Om(H)$ be the Burnside group of $H$ that is, the free abelian group with generators
the various conjugacy classes of subgroups of $H$. To any
finite set $X$ with an $H$-action (or $H$-set)
we can associate an element $(X)\in\Om(H)$ by the requirement
that $(X\sqc X')=(X)+(X')$ for two finite $H$-sets and $(H/H')=H'$ for any subgroup $H'$ of $H$
where $H/H'$ is an $H$-set under left translation.

Let $M(H)$ be the set of all 
pairs $(s,\r)$ where $s\in H$ and $\r$ is an irreducible representation over 
$\CC$ (up to isomorphism) of the centralizer $Z_H(s)$ of $s$ in $H$; the pairs
$(s,\r)$ are taken modulo $H$-conjugacy.
Let $\CC[M(H)]$ be the $\CC$-vector space with basis $M(H)$.

Now let $X$ be a finite $H$-set. For any $(s,\r)\in M(H)$
the fixed point set $X^s$ has an action of $Z_H(s)$ (restriction of the $H$-action on $X$)
hence we can consider the multiplicity $N_{s,\r}$ of $\r$ in the permutation representation of
$Z_H(s)$ on $X^s$. We set $[X]=\sum_{(s,\r)\in M(H)}N_{s,\r}(s,\r)\in\CC[M(H)]$.
Now $(X)\m[X]$ for any finite $H$-set defines a homomorphism

(a) $\Om(H)@>>>\CC[M(H)]$.

\subhead 2.2\endsubhead
We choose a Borel subgroup $B$ of $F^0$ and a maximal torus $T$ of $B$.
Let $F'=\{g\in F;gBg\i=B,gTg\i=T\}$. Then $F'{}^0=T$ and the obvious map
$F'/T\cz_G@>>>F/F^0\cz_G=\bF$ is an isomorphism.
Let $\cb_e^T=\{\fb\in\cb_e;\Ad(t)\fb=\fb\text{ for all }t\in T\}$.
Now $F'$ acts on $\cb_e^T$ by $g:\fb\m\Ad(g)\fb$. 
This action is trivial on $T\cz_G$ hence it induces an action of $F'/T\cz_G=\bF$ on $\cb_e^T$.

Let $s\in\bF$. Let $\cb_e^{T,s}$ be the fixed point set of the action of
$s$ on $\cb_e^T$. Note that $Z_{\bF}(s)$ acts on 
$\cb_e^{T,s}$ as the restriction of the $\bF$-action on $\cb_e^T$. Hence for
any $i$ there is an induced action of $Z_{\bF}(s)$ 
on $H^i(\cb_e^{T,s},\CC)$. We define an element $\ph_e\in\CC[M(\bF)]$ in which the
coefficient of $(s,\r)\in M(\bF)$ is:

(a) $\sum_i(-1)^i(\text{multiplicity of $\r$ in the 
$Z_{\bF}(s)$-module }H^i(\cb_e^{T,s},\CC))$.
\nl
The following is a strengthening of the statement 0.2 that $\BB_e$ is a discretization of $\cb_e$.

\proclaim{Conjecture 2.3} We have $[\BB_e]=\ph_e\in\CC[M(\bF)]$
\endproclaim

\subhead 2.4\endsubhead
Let $W'$ be the affine Weyl group
corresponding to the dual of the adjoint group of $G$. We have $W'\sub W$.
We can find a finite parabolic subgroup $W''$ of $W'$ and a two-sided cell $c''$ of
$W''$ such that $c''\sub c$ (see \cite{L89, 4.8(d)}); moreover, by \cite{L09, 1.5(b2)}, we
can assume that the finite group $\cg_{c''}$ associated to $c''$ in \cite{L87, 3.5} coincides
with $\bF$. Let $\cf_e$ be the set of subgroups of $\bF=\cg_{c''}$ attached in
\cite{L87, 3.8} to the various left cells of $W''$ contained in $c''$ (or
rather one such subgroup in each $\bF=\cg_{c''}$-conjugacy class).
The possible groups $\bF$ are the $\G$ in the list \cite{L20, 0.2(a)}. The subgroups in $\cf_e$
are the subgroups in $\fC(\G)$ with $\G=\bF$, listed in \cite{L20, 1.6}. From \cite{L19} we see
that:

(a) The elements $[\bF/H]\in \CC[M(\bF)]$ for various $H\in\cf_e$ are linearly independent.
\nl
For $b\in\BB_e$ let $\bF_b\sub\bF$ be the stabilizer of $b$ for the $\bF$-action on $\BB_e$.

\proclaim{Conjecture 2.5} $\cf_e$ (see 2.4) is a set of representatives for the
$\bF$-conjugacy classes of subgroups of $\bF$ of the form $\bF_b$ for some $b\in\BB_e$.
\endproclaim

\subhead 2.6\endsubhead
Assuming that 2.3 and 2.5 hold, we see that
the element $(\BB_e)$ of the Burnside group $\Om(\bF)$ is
explicitly determined. Indeed, the element $\ph_e\in\CC[M(\bF)]$ can be explicitly computed
from the knowledge of Green functions for $G$ and its subgroups. Using 2.3 we see that
$[\BB_e]\in\CC[M(\bF)]$ is explicitly determined. Using 2.5 we see
that $(\BB_e)$ is determined by $[\BB_e]$ hence is also explicitly determined.

\subhead 2.7\endsubhead
Assuming 2.5 and that $\r$ is as in 0.1(f) we see that to any $\G\in R(c)$ one can
attach a subgroup $H_\G\in\cf_e$ characterized by the condition that $H_\G$ is conjugate to
$\bF_b$ for some/any $b\in\r\i(\G)$. We note that the subgroups $H_\G\sub\bF$ associated to the
various $\G\in R(c)$ can be regarded as affine analogues of the finite groups associated in
\cite{L87} to the right cells (or left cells) inside a two-sided cell of a finite Weyl group.

\subhead 2.8\endsubhead
For $\x\in\Xi_e$ we denote by $\bF_\x$ the stabilizer of $\x$ in the $\bF$-action on $\Xi_e$.
Assuming 2.5 and the truth of the conjectures in 1.3, we note that the map
$\s:\BB_e@>>>\Xi_e$ in 1.7 is $\bF$-equivariant. Hence if $b\in\BB_e$ then

(a) $\bF_b\sub\bF_{\s(b)}$.
\nl
This seems to be an equality in many (but not all) cases.
Assume for example that $G$ is of type $E_8$ and $e$ is such that $\bF=S_5$.
In this case the subgroups $\{\bF_\x;\x\in\Xi_e\}$ of $\bF$ are exactly the
conjugates of the subgroups in $\cf_e$ (a result of \cite{DLP}); we expect
that in this case (a) is an equality.

Assume now that $G$ is of type $E_8$ and $e$ is of type $E_8(b_6)$ 
(notation as in \cite{CA, p.407}). In this case we have
$\bF=S_3$ and for $\x\in \Xi_e$, $\bF_\x$ is one of the subgroups $S_2,S_3$ or a cyclic group
of order $3$ of $S_3$ (this can be deduced from \cite{DLP, 4.1}); 
if in (a), $\bF_{\s(b)}$ is cyclic of order $3$, we expect to have
$\bF_b=\{1\}$ so that (a) is not an equality.

\head 3. The bijection 0.3(a); an example\endhead
\subhead 3.1\endsubhead
In this section we consider the example where
$G$ is of type $G_2$ and that $e$ is a subregular nilpotent element. 
Let $W$ be as in 0.1. The simple reflections in $W$ are $s_0,s_1,s_2$ where
$s_0s_1$ has order $3$, $s_1s_2$ has order $6$ and $s_0s_2=s_2s_0$.
In this case $c$ is the two-sided cell of $W$ containing $s_0,s_1,s_2$. It is known \cite{L83}
that $c$ consists of all non-identity elements of $W$ with a unique reduced expression.
We write $i_1i_2i_3\do$ instead of $s_{i_1}s_{i_2}s_{i_3}\do$. The elements of $c$ are

$$\matrix 0&01&012&0121&01212&012121&0121210\\
          {}&{}&{}&01210{}&{}&{} \endmatrix$$

$$\matrix 10&1&12&121&1212&12121&121210\\
          {}&{}&{}&1210{}&{}&{} \endmatrix$$ 

$$\matrix 2&21  &212&2121&21212\\
          {}&210&{} &21210&{}\endmatrix$$

Note the apparition of two Coxeter graph of affine type $E_7$
and one of affine type $D_6$.
We write the elements of $\BB_e$ as $[0],[1],[2],[2'],[2'']$.
where the action of $F=\bF=S_3$ on $\BB_e$ keeps $[0]$ and $[1]$
fixed and permutes cyclically $[2],[2'],[2'']$.
The irreducible representations of $S_3$ are denoted by $1,r,\e$ where $r$ is
$2$-dimensional and $\e$ is the sign.
The irreducible representations of $S_2$ are denoted by $1,\e$ where $\e$ is the sign.
The unit repesentation of $S_1$ is denoted by $1$.
The elements of $c$ correspond to the irreducible 
$F$-vector bundles on $\BB_e\T\BB_e$ (up to isomorphism) which appear in the same position in
the following list.

$$\matrix ([0][0];1)&([0][1];1)&([0][2];1)&([0][1];r)&([0][2];\e)&([0][1];\e)&([0][0];\e)\\
          {}&{}&{}&([0][0];r){}&{}&{} \endmatrix$$                                 

$$\matrix ([1][0];1)&([1][1];1)&([1][2];1)&([1][1];r)&([1][2];\e)&([1][1];\e)&([1][0];\e)\\
          {}&{}&{}&([1][0];r){}&{}&{} \endmatrix$$
          
$$\matrix ([2][2];1)&([2][1];1)&([2][2'];1)&([2][1];\e)&([2][2];\e)\\    
           {}&([2][0];1)&([2][0];\e)&{}            \endmatrix$$    

Here a symbol $([?][?];?)$ represents a vector bundle on $\BB_e\T\BB_e$: the first two components
give a point in the support of the vector bundle, the third component is the representation
of the stabilizer of that point in the fibre at that point.

\enddocument